\documentclass[11pt]{article}
\usepackage{amsmath, amssymb, theorem, latexsym, epsfig, amscd}
\numberwithin{equation}{section}

\usepackage{aeguill}
\usepackage{tikz-cd}
\usepackage{environ}
\usetikzlibrary{cd}
\theoremstyle{plain}
\theorembodyfont{\itshape}
\newtheorem{theorem}{Theorem}[section]

\newtheorem{lemma}[theorem]{Lemma}

\theorembodyfont{\rmfamily}
\newtheorem{definition}[theorem]{Definition}

\newtheorem{remark}[theorem]{Remark}

\def\cc{\mathbb C}

\def\wt#1{\widetilde#1}

\def\wh#1{\widehat#1}

\def\kgs {Kleinian groups }
\def\kg {Kleinian group }

\def\ph{\varphi}
\def\te {there exist }
\def\tes {there exists }
\def\st {such that }
\def\qgs {quasi-Fuchsian groups }
\def\qg {quasi-Fuchsian group }
\def\qf {quasi-Fuchsian  }
\def\be {\begin{equation}}
\def\ee {\end{equation}}
\def\cc {\mathbb C}
\def\Om {\Omega }
\def\a {\alpha }
\def\D {\Delta}
\def\d {\delta}
\def\t {\tilde }
\def\homeo {homeomorphism }
\def\g {\gamma}

\def\nbds {neighborhoods }
\def\si {\sigma}
\def\nbds {neighborhoods }

\def\ga {\gamma }

\def\ph {\varphi}

\def\orb{\mbox{orb }}

\title{Uniqueness theorem for completely non-degenerate B-groups}
\author{A.A.Glutsyuk\thanks{CNRS, UMR 5669 (UMPA, ENS de Lyon), Lyon, France}
\thanks{HSE University, Moscow, Russia}
\thanks{Higher School of Modern Mathematics MIPT}
\thanks{The research of A.A.Glutsyuk  is supported by the MSHE project No. FSMG-2024-0048 and by the grant 24-7-1-15-1
 of the Theoretical Physics and Mathematics Advancement Foundation "BASIS".}, Yu.S.Ilyashenko$^\dagger$\thanks{Independent University of Moscow}\thanks{Results of the project "Symmetry. Information. Chaos.", carried out within the framework of the Basic Research Program at HSE University in 2025, are presented in this work.}}
\begin{document}
\newcommand{\bbR}{\mathbb{R}}
\newcommand{\bbP}{\mathbb{P}}
\newcommand{\bbQ}{\mathbb{Q}}
\newcommand{\bbZ}{\mathbb{Z}}
\newcommand{\bbC}{\mathbb{C}}
\newcommand{\bbN}{\mathbb{N}}
\newcommand{\mF}{\mathcal{F}}
\newcommand{\mA}{\mathcal{A}}
\newcommand{\mB}{\mathcal{B}}
\newcommand{\mO}{\mathcal{O}}
\newcommand{\mG}{\mathcal{G}}
\newcommand{\mL}{\mathcal{L}}
\newcommand{\mgF}{\mathgoth{F}}
\newcommand{\pd}{\partial}
\newcommand{\vpf}{\varphi}
\maketitle
\vspace{-1.2cm}
\begin{flushright}
	\textit{In memory of Andrey Andreevich Bolibruch, \\ a brilliant mathematician and a dear friend}
\end{flushright}

\begin{abstract}
  We prove that a completely non-degenerate B-group is uniquely determined by its factor: two such groups with conformally equivalent factors are M\"obius conjugate. A similar property is inherent to the \qgs but not to degenerate B-groups. We also study the factor of a B-group as a triple: the main factor, the marked characteristic complex, and a homotopy class of maps of the first to the second one.
\end{abstract}
\tableofcontents
\section{Introduction}

This paper is devoted to the study of the completely non-degenerate B-groups
introduced  by Bers \cite{Bers1970}
whose definition is recalled below.
They are also called non-degenerate by Bers \cite[pp. 599--600]{Bers1970} and regular by Abikoff \cite{Abikoff1975}.
Regular groups were extensively studied in \cite{Abikoff1975}. We continue the investigations of Maskit and Abikoff.

The factor-spaces of completely non-degenerate B-groups were completely described by Maskit \cite{maskit}. Some properties of these factors were stated but not proved by Maskit. We did not find the proofs of these properties in the literature. So we include it here with a clear understanding that they are known to experts yet not published up to now. Therefore, the paper has a reasonable survey component.

Our first result is the uniqueness theorem. It claims that two B-groups with conformally equivalent factors are M\"obius conjugate. This result is close to \cite{Abikoff1975} and \cite{marden}. Our proof strongly relies upon the Mardens results. On the other hand
the uniqueness theorem was stated in \cite{I}; the proof of this theorem based on the first results of Bers \cite{Bers1970} and Maskit \cite{maskit} on the boundary points of Teihmuller spaces was sketched there. A detailed version of this proof was never published. \footnote{ It is appropriate to give the following historical comment here. The paper \cite{I} ended by a wrong statement: a homeomorphism of the Riemann sphere that is quasiconformal almost everywhere may be quasiconformally extended to the whole spere. This error was mentioned in a review by Abikoff \cite{abrew}. But this error does not break the proof in \cite{I}. It was a wrong justification of a correct statement that some limit map is quasiconformal. In fact, this map is a uniform limit of uniformly quasiconformal maps. Such a limit is quasiconformal itself.}

Another result of our paper is a description of the correspondence between the points of an orbit of a regular B-group and the non-invariant components\footnote{Following B.Maskit \cite{maskit}, by a {\it component} of a Kleinian group $G$ we mean a connected component of its  discontinuity set $\Omega(G)$.} of this group.  In particular, we prove the second statement of theorem 5 of Maskit \cite{maskit}. This correspondence is used in our proof of the uniqueness theorem.

The uniqueness theorem of this paper is crucial for a simultaneous uniformization theorem for algebraic curves with variable topology, subject of a paper in preparation.

\section{Preliminaries}

\subsection{B-groups}

\begin{definition}  \cite[Subsection 1.1.1]{KAG1986} A {\it \kg}  is a group $G$ of M\"obius transformations that acts discontinuously on some nonempty open subset of the Riemann sphere. The maximal open subset where $G$ acts discontinuously is called the {\it discontinuity region} and denoted by $\Om (G)$.
\end{definition}

In what follows, all the \kgs are finitely generated. The complement $\Lambda (G)$ to $\Om (G)$ is called the set of the limit points. It was proved simultaneously and independently by Ian Agol and in the joint work by Danny Calegari and David Gabai \cite{agol, cg}  that the Lebesgue measure of $\Lambda (G)$ is zero. (This was a solution of the famous Ahlfors Measure Conjecture.) For the regular B-groups it was proved much earlier in \cite{Abikoff1975}.

\begin{definition}  A \kg is called a {\it $B$-group} provided that its discontinuity region contains an invariant component $\Delta_0$ homeomorphic to a disc.
We exclude the case of a cyclic group generated by a parabolic transformation with unique fixed point, for which
the invariant component is its complement conformally equivalent to $\cc$. Thus, in our assumptions, $\Delta_0$ is conformally equivalent to the unit disc.
\end{definition}

The factor space $\Om (G)/G$ of a $B$-group is described in \cite{Bers1970} and \cite{maskit}. We reproduce here this description; it will be given for
\emph{torsion free} groups, which  means that $G$ has no elliptic elements.
%
Recall that the factor space
$$
 \Delta_0 /G = S_0
$$
is a Riemann surface of a finite type $(g,n)$: it has genus $g$ and $n$ punctures, by Ahlfors' Finiteness Theorem, see \cite[subsection 1.5.1]{KAG1986}.

 A $B$-group  is called {\it quasi-Fuchsian,} if it has two invariant components. In this case it is known to be quasiconformally conjugated to a Fuchsian group. The class of groups considered here was introduced by Bers and Maskit.

 \begin{definition} \label{def:nondeg}
 Let $G$ be a finitely generated Kleinian group with a single invariant component $\D_0$.
 The group $G$ is called {\it degenerate,} if  $\Omega(G)=\D_0$, i.e., $\Omega(G)$ has no other components; otherwise $G$ is called {\it non-degenerate}.
  \end{definition}
   In other words, a $B$-group is degenerate, if its discontinuity region is simply connected: just one component.
  \begin{remark} The above  terminology, which we use in the paper, is due to B.Maskit \cite{maskit}. Bers \cite[pp.599--600]{Bers1970} uses a different name.
 Namely, let $A$ denote the Poincare area of $\Om(G)/G$, $A_0$ that of $\D_0/G$. It is known that $A \le 2A_0$, and if $G$ is degenerate, then $A = A_0$. Bers calls $G$ non-degenerate if $A = 2A_0$, partially degenerate if $A_0 < A < 2 A_0$. Thus, the groups partially non-degenerate in the sense of Bers are non-degenerate in the sense of Maskit. Abikoff \cite{Abikoff1975} calles the groups non-degenerate in sense of Bers \emph{ regular}, and gives a different definition for them. He proves the equivalence with the Bers' definition; the proof of equivalence is non-trivial. In the present paper we call the latter groups with $A=2A_0$  completely non-degenerate, see Definition \ref{def:nondeg} below.
 \end{remark}

It is known that if a $B$-group is not degenerate and  not quasi-Fuchsian, then it  has a unique invariant component
$\Delta_0$, and it contains the so-called {\it accidental parabolic transformations} (briefly, APTs): those M\"obius transformations that are parabolic
but whose restrictions to $\Delta_0$ are conjugated by uniformization $w: U\to\Delta_0$ to
hyperbolic transformations of the upper half-plane $U$. An APT is called {\it primary,} if it is not a power of
another element of $G$. The set of conjugacy classes in $G$ of primary APTs is known to be finite and is called a {\it basis} of conjugacy classes of APTs. The number of  elements of a basis is called the {\it index} of the group $G$.
See \cite[sections 4, 5]{maskit}.

\subsection{Associated $2$-complex and factors of $B$-groups} \label{sub:ass}

The next  definitions of admissible partition, associated (marked) 2-complex and its factor are introduced in \cite[sections 5 and 7]{maskit}.

\begin{definition}  A set of simple loops $\alpha_1,\dots,\alpha_k$ on a Riemann surface $S$ is called \emph{homotopically independent} or \emph{admissible} provided that:
\begin{enumerate}
\item $\alpha_j$ are pairwise disjoint
\item $\alpha_j$ does not bound a disc or a punctured disc
\item $\alpha_j$ is not freely homotopic to $\alpha_i^{\pm 1}$ for $i \ne j$.
\end{enumerate}
\end{definition}

\begin{definition}  An admissible partition of $S$ is the partition of $S - \{ \a_1, \dots, \a_k\} $ to its connected components $S_1', \dots , S_m'$. Here $\a_1,\dots,\a_k\subset S$ is a homotopically independent
set of loops.\end{definition}

\begin{remark} \label{rbaer} Each admissible set of loops $\alpha_1,\dots,\alpha_k$ and the corresponding admissible partition are uniquely determined by the collection of free homotopy classes of the loops $\alpha_i$ up to isotopy of the ambient surface $S$. This follows from \cite[zusatz 1 p.113, zusatz p.114]{baer} and results of \cite{epstein}.
\end{remark}

For any $j$ let us take an abstract disc $D_j$ bounded by $\a_j$. Then squeeze it to a point. The components $S_l'$ with these squeezed disks added become compact surfaces with punctures and some self intersection points; each $\a_j$ produces one such point. Replace it by a one-dimensional cell $c_j$.  The added closed cells $c_1, \dots , c_k$ are called \emph{connectors}. The resulting $2$-complex is called \emph{an associated complex} and denoted

$$
K = (S_0, \a_1, \dots , \a_k).
$$
The components $S_l$ of the difference
$$
K - \cup \bar c_j = S_1 + \dots + S_m
$$
are called {\it components} of $K$. Some of $S_l$ may be endowed by a complex structure
making it a Riemann surface $S_l^+$ of finite type  so that the canonical identification $I_l:S_l'\to S_l^+$ is orientation-reversing with respect to the orientations defined by the complex structures of the Riemann surfaces $S_0$ and
$S_l$. This additional structure (complex structure collection) is called \emph{marking} of $K$.

Let $G$ be a non-degenerate $B$-group, $S_0=\Delta_0\slash G$. The free homotopy class of each
 loop in $S_0$ defines a conjugacy class in $G$ as in the group of desk transformations of the universal covering $\Delta_0\to S_0$: see the corresponding background material below. It is known that the basis of conjugacy classes of APTs
 (if any) is realized by an admissible (i.e., homotopically independent) collection of loops $\alpha_1,\dots,\alpha_k$, where $k$ is the index of the group $G$ \cite[section 5]{maskit}. We associate to $G$ and
 $\alpha_1,\dots,\alpha_k$
 the above-constructed topological $2$-complex $K$ (at the moment, without any marking).

The following theorem describes factors of the B-groups. The subsequent theorem claims that these factors may be all realized.  Both theorems  are borrowed from \cite{maskit} and adjusted to the torsion free case. The first one is very slightly edited.

To state the first theorem, let us recall the following background material. Let $S$ be an arbitrary Riemann surface, and $\pi:\Delta\to S$ be a universal covering. Fixing base points $o\in S$, $O\in\pi^{-1}(o)\subset\Delta$ defines a natural isomorphism $h_O:\pi_1(S,o)\to G$ between the fundamental group $\pi_1(S,o)$ and the desk transformation group $G$ of the covering as follows. For every $\g \in\pi_1(S,o)$
the corresponding desk transformation $h_O(\g)$  sends $O$ to
 the endpoint of the lifting to $\Delta$ of the loop representing $\g$ with starting point $O$. In the next theorem
 we deal with $\Delta=\Delta_0$, $S=S_0$ and the $B$-group $G$ acting by desk transformations of
 the covering $\Delta_0$.

\begin{theorem}   \label{thm:fact}  Let $G$ be a $B$-group with invariant component
$\D_0$. Let $S_0 = \D_0 / G$ and let $K = K(S_0, \a_1,
\dots , \a_k)$ be the associated $2$-complex. Then there is a marking $S_1^+, \dots , S^+_s$ on $K$, so that $\Om (G)/ G = S_0 + S_1^+ + \dots + S^+_s$. Furthermore, if $\D $ is any  component of $G$ different from $\Delta_0$, then there are points $o$ on $S_0$ that does not lie on any $\a_i$, and $O$ on $\D_0$ projected to $o$,  so that $G_\Delta$ (the stationary subgroup of $\D $ in $G$) is the image, under the natural homomorphism $h_{O}:\pi_1(S_0,o)\to G$, of the subgroup of $\pi_1(S_0, o)$ generated by loops which do not cross any $\a_i$.
\end{theorem}

 This is theorem 5 from \cite{maskit}. Surprisingly, only its first statement  is proved in \cite{maskit}. There is no one word in \cite{maskit} of the proof of the second statement. We did not succeed to find the proof of this statement in literature. So we include it here. In fact, this proof is written in \cite{maskit} between the lines; we need only to read it and write down. Here and there we use deliberately and literally some lines from \cite{maskit}.

 \begin{definition}\label{def:nondeg}
   If all the components $S_1,...,S_m$ of the 2-complex mentioned in Theorem \ref{thm:fact} are marked, then the corresponding group is {\it completely non-degenerate.}
 \end{definition}

 The next theorem, which is theorem 6 from \cite{maskit}, is in a sense inverse to the previous one. We include it here for the sake of completeness and for future references.

 \begin{theorem}   \label{thm:real} Let $S_0$ be a finite Riemann surface of type $(g_0, n_0)$, where $3g_0 - 3 + n_0 >0$. Let $\a_1, \dots , \a_k$ be a homotopically independent set of loops on $S_0$, and let $K = K(S_0, \a_1, \dots , \a_k)$ be the associated $2$-complex.

Let $S_1^+, \dots S^+_s$ be some marking on $K$. Then there is a $B$-group $G$ with invariant component $\D_0$ so that
\begin{enumerate}
\item $\D_0 / G = S_0$,
\item $\alpha_1,\dots,\alpha_k$ correspond to the basis of conjugacy classes of primary APTs in $G$,
\item $\Om (G) / G = S_0 + S^+_1 + \dots + S^+_s$.
\end{enumerate}
\end{theorem}

 \subsection{Factors of \qgs and the uniqueness theorem for such groups} \label{sub:factqf}

 The well known material of this subsection (due to Bers \cite{Bers1960}) serves as a prototype of our main result.

 By definition, a \qg is a B-group whose discontinuity set consists of exactly two invariant components; denote them by $\D_0$ and $\D_1$. Its factor-space is
 $$
 \Om(G)/G = S_0 + S_1,
 $$
 where $S_0= \D_0/G $ and $S_1 = \D_1/G $ are homeomorphic Riemann surfaces. Let $\pi$ be the corresponding projection. Furthermore, a homotopy class $[\ph]$ of orientation-reversing homeomorphisms $S_0 \to S_1$ is well defined, \st roughly speaking, for any $\Phi \in [\ph]$ the natural homomorphisms of $\pi_1(S_0)$ and $\pi_1(S_1) = \pi_1(\Phi(S_0))$ to $G$ agree with each other. In more detail, \te $o \in S_0$, $O \in \D_0, \pi(O) = o$, \st for any $\Phi \in [\ph]$ \tes a $Q \in \D_1$ \st for any $\ga \in \pi_1(S_0,o)$ the following elements of the group $G$ coincide:

\begin{equation}\label{eqn:natur}
 h_O(\ga) = h_Q(\Phi (\ga)).
\end{equation}

 \begin{definition}
 The factor of the \qg $G$ is the triple $(S_0, S_1, [\ph])$ described above.
 \end{definition}

 \begin{definition}
   Two factors $(S_0, S_1, [\ph])$ and $(\t S_0, \t S_1, [\t \ph])$ are conformally equivalent provided that \te conformal maps $\Psi_0: S_0 \to \t S_0, \Psi: S_1 \to \t S_1 $
     \st $ \Psi \circ [\ph] = [\t \ph] \circ \Psi_0 $.
 \end{definition}

 \begin{theorem}
   (uniqueness theorem for \qgs) Two \qgs with conformally equivalent factors are M\"obius conjugate.
 \end{theorem}
 As was mentioned above, these definitiions and the theorem are borrowed from \cite{Bers1960}. Bers calls the factor $(S_0, S_1, [\ph])$ the connected pair.

 \subsection{Factors of  non-degenerate B-groups  and  uniqueness theorem } \label{sub:factnon}

 We pass to a sort of addendum to the Maskit theorem \ref{thm:fact}. The latter theorem  describes completely the factor-space of a B-group, but  some extra information is hidden in this space. The situation of non-quasi-Fuchsian B-group is different and we need to change the definition of factor.

To define the factor of a non-quasi-Fuchsian B-group, we use the following topological construction. 
 Let $G$ be a non-degenerate B-group. Any component of the associated complex corresponding to $G$ may be constructed in the following way. Let us cut the surface $S_0$ along each loop $\alpha_j$. This transforms the complement $S_0\setminus\sqcup_j\alpha_j$ to the interior of a two-dimensional surface with boundary,
denoted $\wh S_0$. Each component $S_\ell'$ of the complement  $S_0\setminus\sqcup_j\alpha_j$ is identified with a connected component of $Int(\wh S_0)$. Each $\alpha_j$ generates  two boundary components of
the cut surface $\wh S_0$: two copies of $\alpha_j$.  Let us attach to them punctured discs denoted by $C_j'$ and $C_j''$: one disk per copy. Let us attach to each $S_\ell'$ all those of the latter punctured discs that are adjacent to $S_\ell'$. The union thus obtained is homeomorphic to a component $S_\ell$ of the associated complex,  and it will be identified with $S_\ell$. For every $j$ fix a narrow  tubular neighborhood $U_j$ of the loop $\alpha_j$ in the initial surface $S_0$. In the new surface $\wh S_0$ with the above disks attached,  the annulus $U_j$ cut along the loop $\alpha_j$  becomes a union of two topological annuli, one adjacent to the disk $C_j'$ and the other to  $C_j''$. The disks $C_j'$, $C_j''$ enlarged by  their latter  adjacent annuli will be denoted by $C^1_j$ and $C^2_j$ respectively.
Set
 $$S_\ell''  = S_\ell' - \cup (C_j^1 + C_j^2),$$
  the union is taken over those $j$ for which $\a_j \subset \partial S_\ell'$. Let $U = \cup U_j$.

 \begin{lemma} \label{lem0}
 In the notations above, \tes a continuous map $\Phi: S_0 \to K$ \st the following holds.

 \begin{itemize}
   \item \begin{equation} \Phi (S_\ell'') \subset S_\ell \ \ \forall \ell \in \{1,...,m \},\label{phiu0}
   \end{equation}
   \item \begin{equation} \Phi (U_j) = C_j^1 + C_j^2 + c_j \ \ \text{ for every } \ j=1,\dots,k,\label{phiu}\end{equation}
   where $c_j$ is the connector corresponding to  $\a_j$.
   \item For every marked component $S_t^+$ of the associated complex the following statements hold:

   \begin{equation}\text{ the map } \ \ \Phi:S_t''\to S_t^+ \ \text{ reverses the orientation}
   \label{reverse}
   \end{equation}
   induced by complex structures
   of the Riemann surfaces $S_t''\subset S_0$, $S_t^+$;

-  for any $o \in S_t''$, $O \in \pi^{-1}(o)\subset\Delta_0$ and $q = \Phi (o)\in S_t^+$ \tes a $Q \in \pi^{-1} (q)$ \st for any $\g \in \pi_1 (S_t'',o)$ the following two elements of $G$ coincide:
       $$  h_O(\g) = h_Q (\Phi (\g))             .$$
 \end{itemize}
 \end{lemma}

 This lemma immediately follows from the Main Lemma stated below.

 \begin{definition}\label{def:factor}
   The {\it factor} of a completely non-degenerate B-group with a characteristic complex
   $K$ with marking is the triple
   \begin{equation}\label{eqn:facto}
     (S_0,  K, [\ph])
   \end{equation}
   where $[\ph] : S_0 \to K$
   is the homotopy class of the map $\Phi$ from Lemma \ref{lem0}.
\end{definition}

\begin{definition}\label{def:confequiv}
   Two factors of  completely non-degenerate B-groups with  marked characteristic complexes $K$ and $\t K$, \eqref{eqn:facto} and
   \begin{equation}\label{eqn:factor}
     (\t S_0, \t K, [\t \ph])
   \end{equation}
   are {\it conformally equivalent} provided that
   \te a biholomorphic map $\Psi_0: S_0 \to \t S_0$ and a \homeo $\Psi: K \to \t K$ \st
   $\Psi (S_t^+) = \t S_t^+$, the diagram
      $$
\begin{CD}
S_0 @> {[\ph ]}>> K\\
@V\Psi_0VV      @VV\Psi V\\
\t S_0 @> {[\t \ph ]}>> \t K\end{CD}
$$
commutes, and the restrictions
 $$
 \Psi_t = \Psi|_{S_t^+}: S_t^+ \to \t S_t^+
 $$
 are biholomorphic.
 \end{definition}

\begin{theorem} \label{thm:uniqu}  Two completely non-degenerate $B$-groups with conformally equivalent factors are M\"obius conjugate.
\end{theorem}

This is the first main result of our paper.

\begin{remark} Quasiconformal deformations and Teichm\"uller spaces of Kleinian groups were studied in papers \cite{kra, maskit2}. Their results imply that the Teichm\"uller space of a completely non-degenerate $B$-group is isomorphic to the product of Teichm\"uller spaces of  stabilizers of  a finite collection of its components that are bijectively projected
to the collection of the components of the quotient. These results themselves do not imply our Theorem \ref{thm:uniqu}. 
\end{remark}

\subsection{Marden's results and Conjugacy theorem}

The proof of Theorem \ref{thm:uniqu}  is based on the two following results by A.Marden (we quote them with the notations slightly changed).

\begin{theorem}  \label{tmarden} (Marden, \cite{marden}, p. 429, Theorem 8.1). Let $G$ and $\wt G$
be Kleinian groups such that

(i) the natural action of the group $G$ on the hyperbolic 3-space has a finite-sided fundamental polyhedron
(in this case the group $G$ is called {\bf geometrically finite});

(ii) there exists a conformal isomorphism $R:\Omega(G)\to\Omega(\t G)$ such that the map
$g\mapsto R\circ g\circ R^{-1}$ yields an isomorphism $H:G\to \t G$ conjugating the action of $G$
on $\Omega(G)$ with the action of the group $\t G$ on $\Omega(\t G)$.

\noindent Then $H$ is an inner automorphism: the groups $G$ and $\t G$ are conjugated by a M\"obius transformation.
\end{theorem}

\begin{remark} Theorem \ref{tmarden} is stated in \cite{marden} in a more general form, under a weaker condition, when
$R$ is just an orientation preserving homeomorphism: it states that then the groups are quasiconformally
conjugated, and if $R$ is conformal, then $H$ is an inner automorphism.
\end{remark}

\begin{theorem} \label{tmarden2} (Marden, \cite{marden}). For every completely non-degenerate $B$-group
its natural action on the hyperbolic 3-space has a finite-sided fundamental polyhedron.
\end{theorem}

Theorem \ref{tmarden2} follows from \cite[proposition 4.2]{marden}.

Theorem \ref{thm:uniqu} follows from Theorems \ref{tmarden}, \ref{tmarden2} and the next theorem proved below.

\begin{theorem} \label{conjd} (Conjugacy theorem) Two completely non-degenerate $B$-groups with equivalent factors are conformally conjugate on their discontinuity sets.
\end{theorem}

First we prove Theorem \ref{conjd} modulo the Main Lemma stated below. Afterwards we prove the Main Lemma.

\section{The Main Lemma on non-degenerate $B$-groups}

Recall that  we consider the groups without the elliptic elements. So the superscript $'$ used in \cite{maskit} is omitted.

The Main Lemma stated below is a refinement of the second statement of theorem 5 in \cite{maskit}, see Theorem \ref{thm:fact} above.

  \begin{lemma}
   (Main Lemma) \label{lem:main}
   Let $G$ be a non-degenerate $B$-group with the invariant component $\D_0$,
   $K$ be the associated marked complex. Let other notations from Section \ref{sub:factnon} hold. Then there exists a homotopy class of
  continuous maps $[\ph]: S_0 \to K$ with a representative $\Phi$ satisfying (\ref{phiu0}),
   (\ref{phiu}), (\ref{reverse}) (a representative satisfying (\ref{phiu0}), (\ref{phiu}), (\ref{reverse}) will be called {\bf good}) such that each  good representative $\Phi$  satisfies the following statements for every $t$ corresponding to a marked component $S_t^+$.
   (In the case, when $G$ is quasi-Fuchsian, we claim that the corresponding
   projection $S_0\to S_1$ introduced in Subsection 2.3 satisfies these statements.)

   A) For every $o\in S_t'$ for any component $\D \subset \pi^{-1} (S_t^+)\subset\Omega(G)$ there exists a point $O \in \D_0$ \st $\pi (O) = o$, and
   \begin{equation}\label{eqn:stab}
    h_O (\pi_1 (S_t',o)) = G_\D:=\text{ the stabilizer of the component } \Delta  \text{ in } G;
  \end{equation}
 and for any $Q \in \D, q = \pi (Q) \in S_t^+$ one has
   \begin{equation}\label{eqn:stab1}
     h_Q (\pi_1 (S_t^+,q) ) = G_\D.
  \end{equation}
   For every  good representative $\Phi\in [\ph]$ with\footnote{Here we consider that the tubular neighborhoods $U_j$ around loops $\a_j$ defining $\Phi$ are small enough so that $o\in S_t''$, thus $\Phi(o)$ is well-defined and lies in $S_t^+$, by (\ref{phiu0}).}
   $q = \Phi (o)\in S_t^+$

   \tes a $Q \in \D, \pi (Q) = q$
   \st for any $\g \in \pi_1 (S_t'', o)=\pi_1(S_t',o)$
   one has
   \begin{equation}\label{eqn:pione}
    h_O(\g)  = h_Q(\Phi(\g)).
  \end{equation}

  B) Let $O$ and $Q$ be the same as in A. Then there is a map $\si: \orb_G O \to \orb_G Q$, $\sigma(O)=Q$, such that for any $O' \in \orb_G O$ one has
  \begin{equation}\label{eqn:iso}
    h_{O'}(\pi_1(S_t',o)) =  h_{\si(O')}(\pi_1(S_t^+,q)) = G_{\D'}
  \end{equation}
  Here $\D'$ is the component of $G$ that contains $\si(O')$. Moreover,
identity (\ref{eqn:pione}) holds with $O$, $Q$, $\Delta$ replaced by $O'$, $\si(O')$ and
  $\Delta'$ respectively.
 The map $\si$ commutes with the $G$ action:
  the diagram

  \begin{equation}\label{eqn:upper}
\begin{CD}
 \orb_G O @>G>>  \orb_G O\\
 @VV{\sigma}V            @V{\sigma}VV\\
 \orb_{ G}  Q  @>G>>  \orb_{G} Q                                               \end{CD}
\end{equation}

  is commutative.
  \end{lemma}

  This is the second main result of our paper.

  \begin{remark}
    The lemma implies that the map $\si^*: O'\to\Delta'$,
    $\sigma^*:\orb_G O \to \{\mbox{components of } \pi^{-1} (S_t^+)\}$
    is an epimorphism.
  \end{remark}

  \begin{remark}
    Statement B is a trivial consequence of A, and is given for the future convenience.
  \end{remark}

\begin{remark}
Lemma \ref{lem0} is a direct consequence of the main one.
\end{remark}

\section{Proof of Theorem \ref{conjd} modulo the Main Lemma}
This section is devoted to the proof of Theorem \ref{conjd} and hence, the uniqueness theorem  modulo the Main Lemma.

\subsection{Definition of the conjugacy in the main component}

In what follows, we will use many times the lift of the Riemann surface to its universal cover. Let $S$ be a Riemann surface, and $\D$ be its universal cover with the projection $\pi: \D \to S$. The inverse map $\pi^{-1}$ is not at all uniquely defined; two such maps differ by a desk transformation. But if we choose an arbitrary point $o \in S$ and a point $O \in \pi^{-1}(o)$, then the lift of $\pi$ that brings $o$ to $O$ is uniquely defined; this is the extension of the germ $\pi^{-1}: (S,o) \to (\D, O)$; we denote this extension by $\pi^{-1}_{O,o}$.

Let us turn to the proof of the Conjugacy Theorem \ref{conjd}, namely to the construction of the conformal map $R: \Om(G) \to \Om(\t G)$. Let us first construct it on $\Delta_0$.  Let, as before,
$\pi$ be the projection $\Om( G) \to \Om( G)/  G$ and define $\t \pi: \Om(\t G) \to \Om(\t G)/ \t G$.
Take an arbitrary point $O \in \D_0, o = \pi(O), \t o = \Psi_0 (o)$, and take an arbitrary point
$\t O \in \t \pi^{-1}(\t o)$. Let us define $R|_{\D_0}$ by
\begin{equation}\label{eqn:rdes}
  R = \t \pi^{-1}_{\t O, \t o} \circ \Psi_0 \circ \pi.
\end{equation}
Obviously, the map $R$ conjugates the action of $G$ on $\D_0$ with the action of $\t G$ on $\t \D_0$. Define
\begin{equation}\label{eqn:hconj}
  H: G \to \t G, \ \ g \mapsto H(g) = R \circ g \circ R^{-1}.
\end{equation}
In fact, to define the element $H(g)$ it is sufficient to define its action on the orbit of one  point $O_t$ (even on one point only because we know that $H(g)
$ belongs to $\t G$). So $H: G \to \t G$ is defined by the commutative diagram

\begin{equation}\label{eqn:central}
\begin{CD}
 \orb_G O_t @>G>>  \orb_G O_t\\
 @VVRV            @VVRV\\
 \orb_{\t G} \t O_t  @>H(G)>>  \orb_{\t G} \t O_t                                               \end{CD}
\end{equation}
Here $\wt O_t=R(O_t)$.

\subsection{ Definition of the conjugacy in the complement $\Om \setminus \D_0$   }

It is sufficient to define $R$ on $  \pi^{-1}(S_t^+)$ for any  $ t \in \{1,...,s \}$ as a lift of the conformal map $\Psi_t:S_t^+\to\wt S_t^+$. Let $\Phi:S_0\to K$ be a map  from the Main Lemma.  Take an arbitrary $o_t \in S_t'$. We can and will consider that $o_t\in S_t''$,
thus $q_t:=\Phi(o_t)\in S_t^+$, see Footnote 3.   Take a point $O_t \in \pi^{-1}(o_t)$. Let $Q_t$  and the map $\si: \orb_G O_t \to  \orb_G Q_t $ be the same as in the Main Lemma. Take the similar objects for the group $\t G $ denoted in the same way, but with tilde, with $\wt o_t=\Psi_0(o_t)$ and $\wt\Phi=\Psi\circ\Phi\circ\Psi_0^{-1}$,
$\Psi=(\Psi_1,\dots,\Psi_s)$: see the definition of conformally equivalent factors.
Then one has $\wt q_t=\wt\Phi(\wt o_t)$.

For any component $\D_t \subset \pi^{-1}(S_t^+)$ we have to find a component
$\t \D_t \subset \t \pi^{-1}(\t S_t^+)$ that will be the target of $R|_{\D_t}$.
This is done with the use of the Main Lemma.

Take any point $O' \in \orb_G O_t$ \st $\si^* (O') = \D_t$. This point exists but is not at all unique. The point $O'$ belongs to $\D_0$ where $R$ is already defined. Take $\t O' = R(O')$ and $\t Q' = \t \si (\t O'), \ \t \D_t = \si^* (\t O')\owns \t Q'$. An equivalent definition of $\t Q'$ may be seen from the following commutative diagram:

\begin{equation}\label{eqn:left}
\begin{CD}
 \orb_G O_t @>\si>>  \orb_G Q_t\\
 @VVRV            @VRVV\\
 \orb_{\t G} \t O_t  @>\t \si>>  \orb_{\t G} \t Q_t                                               \end{CD}
\end{equation}

The right arrow is the definition of $R$ on
$\orb_G Q_t.$

Let us now take
\begin{equation}\label{eqn:rdet}
  R|_{\D_t} = \t \pi^{-1}_{\t \si (\t O'), \t q_t} \circ \Psi_t \circ \pi.
\end{equation}

The definition depends on the arbitrary choice of $O' \in \pi^{-1}(o_t)$. We have to prove that $R$ is well defined (does not depend on this choice), and conjugates the action of $G$ and $\t G$; namely, the diagram
\begin{equation}\label{eqn:leftt}
\begin{CD}
 \Om - \D_0 @>G>>  \Om - \D_0\\
 @VVRV            @VRVV\\
 \t \Om - \t \D_0  @>H(G)>>  \t \Om - \t \D_0
  \end{CD}
\end{equation}
commutes. As the lifting of a Riemann surface to its universal cover is uniquely defined by the choice of the base point, it is sufficient to prove the commutativity of the following diagram:

\begin{equation}\label{eqn:large}
\begin{CD}
 \orb_G Q_t @>G>>  \orb_G Q_t\\
 @VVRV               @VRVV\\
 \orb_{\t G} \t Q_t  @>H(G)>>  \orb_{\t G} \t Q_t
  \end{CD}
\end{equation}

\subsection{ Conjugacy relation. Proof of Theorem \ref{conjd}}

To prove the commutativity of this diagram, consider a large diagram shown below.

\bigskip

 \hskip1.5cm\begin{tikzcd}
 \orb_GQ_t   \arrow[rrr, "G"]  \arrow[ddd, "R" '] & & & \orb_GQ_t \arrow[ddd, "R"] \\
 & \arrow[ul, "\sigma" '] \orb_GO_t  \arrow[r,"G"] \arrow[d, "R" '] &  \orb_GO_t \arrow[ur, "\sigma"]   \arrow[d, "R"]\\
 & \arrow[ld, "\wt\sigma"] \orb_{\wt G}\wt O_t \arrow[r, "H(G)" '] & \orb_{\wt G}\wt O_t
 \arrow[dr, "\wt\sigma" '] \\
 \orb_{\wt G}\wt Q_t \arrow[rrr, "H(G)" '] & & & \orb_{\wt G}\wt Q_t
 \end{tikzcd}

\bigskip

This diagram is a union of five small commutative diagrams.
The one in the center is the diagram \eqref{eqn:central}.
The one on the left is \eqref{eqn:left}.
The one on the right is the same as the left one.
The upper one is \eqref{eqn:upper}. The lower one is the same as the upper, but for the group $\t G = H(G)$ instead of $G$. All these diagrams commute. This implies the commutativity of all the large diagram. In particular, the small diagram formed by the arrows on the edges forms the diagram \eqref{eqn:large}. This proves its commutativity as required. This commutativity in turn implies that the map $R$ is well defined. This proves Theorem \ref{conjd} modulo the Main Lemma.

\section{Proof of the Main Lemma}

\subsection{  Statement B}
Let $\D, O$ and $Q$ be the same as in A, and let \eqref{eqn:stab} hold.
 Let $O' = g (O) \in \orb_G O, \ g \in G, \ Q = \sigma (O)$. Then let $\si (O') = g (Q)$.

Let $\D' = g(\D).$ Then

$$G_{\D'} = g \circ G_\D \circ g^{-1}.$$
On the other hand,

$$
h_{g(O)} (\pi_1 (S_t',o)) = g \circ h_O (\pi_1 (S_t',o))\circ g^{-1} = g \circ G_\D \circ g^{-1}.
$$
This proves \eqref{eqn:iso}. The diagram \eqref{eqn:upper} is commutative  by construction. This proves B.

Proof of A. Statement \eqref{eqn:stab} is very close to the Maskit theorem \ref{thm:fact}. This statement is proved below. Statement \eqref{eqn:stab1} is a trivial consequence of the definition of lifting. Statement \eqref{eqn:pione} will be proved together with \eqref{eqn:stab}. We prove statement \eqref{eqn:stab} by induction in the number of loops $\a_1,...,\a_k$. At the same time the existence of the class $[\ph]$ with the property \eqref{eqn:pione} is proved.

\subsection{Base of induction: $k = 0 $.}

If $k = 0 $, then the group $G$ has no APT's. Hence, it is quasi-Fuchsian, by theorem 4 from \cite{maskit}. Consider first the case when it is Fuchsian.

For any Fuchsian group $G$, we have: $\Om (G) = L + U$, where $L$ and $U$ are the lower and upper halfplanes respectively. Let $\hat \Phi_0: L \to U$ be the symmetry $z \mapsto \bar z$. Denote: $L = \D_0, \ U = \D_1, \ S_0 = L/G, \ S_1^+ = U/G$. Let $\pi: L + U \to S_0 + S_1^+$ be the natural projection. Take an arbitrary point
$o \in S_0, O \in \pi^{-1}(o), Q = \hat \Phi_0(O), q = \pi (Q)$. Then
$$
G = h_O(\pi_1 (S_0,o)) = h_Q(\pi_1 (S_1^+,q)) = G_{\D_1}
$$
as required.

Let now $ \Phi_0: S_0 \to S_1^+$ be the descending of the symmetry $\hat \Phi_0$:
$\ \Phi_0\circ \pi = \pi \circ \hat \Phi_0$. Let $o, O, Q$ be the same as above.
Then for any $\g \in \pi_1 (S_0,o)$
$$
h_O(\g) = h_Q (\Phi_0 (\g)).
$$

Let now $\Phi$ be a \homeo $S_0 \to S_1^+$ homotopic to $\Phi_0$, set $q'=\Phi(o)$.
We will find $Q'$ \st \eqref{eqn:pione} holds for $Q$ and $\Phi_0$ replaced by $Q'$ and $\Phi$:

\begin{equation}\label{eqn:pionne}
h_O(\g) = h_{Q'} (\Phi (\g)).
\end{equation}

Let $\{\Phi_\tau \ |\ \tau \in [0,1], \Phi_1 = \Phi \}$ be a homotopy from $\Phi_0$ to $\Phi$. Take the  path
$\g_\Phi: [0,1] \to S_1^+$, $\g_\Phi(\tau)=\Phi_\tau(o)$. Let $ \hat \g_\Phi$ be the lift of this path  to $U$ with the starting point $Q$.
Let $Q'$ be the endpoint of this lift. Then  for any $\g \in \pi_1 (S_0,o)$ the relation \eqref{eqn:pionne}
holds as required.

The case of quasi-Fuchsian group is reduced to the case of Fuchsian group, since each quasi-Fuchsian group
is conjugated to a Fuchsian group by a quasiconformal homeomorphism. This proves the induction base.

\subsection{Induction step: case when every loop $\a_j$ divides $S_0$
}

We begin with what is called case 2 in \cite{maskit}, section 10: every loop $\a_j$ divides $S_0$. We repeat the construction from \cite{maskit}, section 10, case 2.

We re-order the $\a_i$, if necessary, so that $\a_k$ divides $S_0$ into two surfaces  $S_{01}$ and $S_{02}$ where every $\a_i, i < k$, lies on $S_{01}$.

Pick a base point $o \in S_0$, where $o \in \a_k$. Let $O \in \pi^{-1}(o)$  be some base point in $\D_0$. For $i = 1, 2$, let $\pi_i^*$ be the subgroup of $\pi_1(S_0,o)$ generated by the loops which lie, except for $o$, in $S_{0i}$. For $i = 1, 2$, let $G_i$ be the image, under the natural homomorphism, of $\pi_i^*$:
 $$
 G_i = h_O(\pi_i^*).
 $$
 For $i = 1, 2$, let $\t \D_{0i}$ be that  connected component of the pre-image, under the natural projection, of $S_{0i}$ which has $O$ on its boundary. For $i = 1, 2$, let $ \D_{0i}^*$ be that component of $G_i$ which contains $\t \D_{0i}$. For each APT $\d \in G$ that corresponds to $\a_k$,  there is a simple arc
 $B_\d' \subset\Delta_0$, where $B_\d'$ is invariant under $\d$ and $B_\d'$ projects onto $\a_k$. Adjoining the fixed point of $\d$ to $B_\d'$, we get a simple closed curve $B_\d$.

Any topological circle $B_\d$ bounds two topological discs $C_{\d 1}$ and $C_{\d 2}$ whose closures together constitute the whole Riemann sphere. Let $C_{\d i}$ be that one of these two discs that does not intersect $\wt\D_{0i}$. Denote the set of all such discs for all $\d$ above by $\mathcal C_{ i}$. In \cite{maskit},  section 10, case 2, it is proved that
$$
\D_{0i}^* = \wt\D_{0i} + \sum_{\mathcal C_{ i}}C_{\d i}+\sum_\delta B_\delta',
$$
  and
  $$
  S_{0i}^*:= \D_{0i}^*/G_i = S_{0i} + C_i,
  $$
where $C_i$ is a punctured disc with boundary, pasted to  $S_{0i}$ along $\a_k$.
Let
$$
\Om(G_1)/G_1 = S_{01}^* + S_1^+ +...+S_{s}^+
$$
Let $K$, $K_1$ denote the marked complexes associated to the groups $G$ and
$G_1$ respectively. Let $\Phi^*:S_{01}^*\to K_1$ be a good homeomorphism satisfying Statement A of the Main Lemma
for the group $G_1$.
 We will show that the
factor $(\Omega(G)\setminus\Delta_0)\slash G$ consists of the surfaces $S_1^+,\dots,S_s^+$ and maybe of an extra surface $S_{s+1}^+$. We construct the corresponding map $\Phi_{s+1}:S_{s+1}''\to S_{s+1}^+$ for the surface $S_{s+1}^+$ and glue it
with $\Phi^*$ to one map $\Phi:S_0\to K$. We prove
the statement of the Main Lemma for the group $G$ and thus constructed map $\Phi$.
As all the curves $\a_1,...,\a_{k-1}$ belong to $S_{01}$, the group $G_2$ is a B-group with no APT's. By Theorem 4 from \cite{maskit}, the group $G_2$ is either quasi-Fuchsian, or degenerate.

In both cases the associated complex of the group $G$ is the union of the associated complex of the group $G_1$, the surface $S_{02}$ with its boundary arc $\alpha_k$ squeezed to a point and the  connector $c_k$ corresponding to $\alpha_k$.
The surface $S_{02}$ is marked (arises in the factor) if and only if
 $G_2$ is quasi-Fuchsian.

Case 1: the group $G_2$ is quasi-Fuchsian. Thus, it has two components: $S_{02}^*$ and another one, denoted by
$S^+_{s+1}$. Moreover, \tes a homotopy class of maps $[\ph_{s+1}]: S_{02}^* \to S^+_{s+1}$ described in Section \ref{sub:factqf}. For any B-group $G$ denote by $\pi_G$ the projection of the discontinuity set of $G$ to its factor; we omit the subscript $G$ when it does not produce misunderstanding. For any $o \in S^*_{02}, O \in \pi_{G_2}^{-1}(o)\subset\Delta_{02}^*$ and any
$\Phi_{s+1}\in [\ph_{s+1}]$ \tes a $Q \in S^+_{s+1}$ \st for any $\g \in \pi_1 (S^*_{02},o)$
\begin{equation}\label{eqn"pionnne}
  h_O(\g) = h_Q (\Phi_{s+1} (\g)).
\end{equation}
This is a property of \qgs that we merely refer to.

Let us now define $[\ph]: S_0 \to K$ in the following way. Take two nested tubular \nbds of
$\a_k$:  $U_k$ and $u_k \Subset U_k$. Set $U_{ik}:=(U_k - u_k)\cap S_{0i}$, $u_{ik}:=u_k \cap S_{0i}$,
$i = 1,2.$ Let $c=c_k$ denote the connector in the factor complex of the group $G$  corresponding to $\alpha_k$. Fix a map $\Phi_{s+1}\in [\ph_{s+1}]$
defined above. Let us now define a map $\Phi$ in the following way:
$$
\Phi =
\begin{cases}
\Phi^* \mbox{ on } S_{01} - U_k\\
\Phi_{s+1} \mbox{ on } S_{02} - U_k\\
\mbox{ a \homeo} U_{ik} \to C_i \mbox{ on } U_{ik}, i = 1, 2 \\
\mbox{ a map } a: u_k \to c \mbox{ defined below}.
\end{cases}
$$

The map $a$ of $u_k$ onto the connector $c$ is defined as follows. Let $p: \a_k \times c \to c$
be a projection along the first factor, and $b: u_k \to  \a_k \times c $ a homeomorphism. Then $a = p \circ b$. Now define $[\ph]$ as the homotopy class of $\Phi$.
This completes the construction of $[\ph]$ in the case, when $G_2$ is quasi-Fuchsian.

Case 2: the group $G_2$ is degenerate. Then $\Omega(G_2)=S_{02}^*$. But  there is no additional component of $\Omega(G)\slash G$ coming from $G_2$: the surface $S_{s+1}$ is not marked. In this case we can take $\Phi_{s+1}:S_{s+1}''\to S_{s+1}$ to be the canonical inclusion given by the construction at the beginning of Subsection 2.4 and then
extend $\Phi_{s+1}$ and $\Phi^*$ to a map $\Phi:S_0\to K$ as above.

Let us now prove \eqref{eqn:stab} and \eqref{eqn:pione}. The proof presented below works in both cases, when $G_2$ is quasi-Fuchsian or degenerate. In \cite{maskit}, Section 10, it was proved that  \te fundamental sets $D$, $D_1$, $D_2$ for the groups $G, G_1, G_2$
respectively in their discontinuity sets \st
\begin{equation}
D -\D_0 = (D_1 - \D_{01}^* ) + (D_2 - \D_{02}^* ).\label{d12}
\end{equation}
This implies that each non-invariant component $\Delta$ of the group $G_i$ different from $\D_{0i}^* $ is  at the same time a non-invariant component of $G$.  Indeed,  every point $x\in\Delta$ is $G_i$-equivalent to a unique point
$y\in D_i\ - \D_{0i}^*$, $y=g(x)$, $g\in G_i$, and hence, $y\in D -\D_0\subset\Omega(G)$, by (\ref{d12}).
This implies that $\Delta\subset\Omega(G)$. The component in $\Omega(G)$ containing
$\Delta$ coincides with $\Delta$, since $\Omega(G)\subset\Omega(G_i)$ by the
inclusion  $G_i\subset G$.

Conversely, each non-invariant component $\Delta$ of $G$ is $G$-equivalent to some non-invariant component of some of the groups $G_1, G_2$.
Indeed, each its point is sent to a point  $y\in D -\D_0$ by some (and unique) element
 $g\in G$. Thus, $y$ lies in
 $D_i-\D_{0i}^*$ for some $i=1,2$ by (\ref{d12}). The component
  of the group $G_i$ containing $y$ is a component of the group $G$ by the above argument, and hence, it coincides with $g(\Delta)$.

 For every non-invariant component $\Delta$ of the group $G_i$ the stabilizer $G_\Delta$ coincides with
its stabilizer  in the group $G_i$. Indeed, two points  $x_1,x_2\in \Delta$ are $G$-equivalent, if and only if the corresponding points $y_j\in D - \D_0$ coincide.
   But there exist  $g_j\in G_i$ such that
   $g_j(x_j)\in D_i- \Delta_{i0}^*\subset D - \Delta_0$. One has
   $y_j=g_j(x_j)$ by uniqueness. Therefore, $x_1$ and $x_2$ are $G$-equivalent, if and
   only if  $x_2=g_2^{-1}g_1(x_1)$, in which case they are $G_i$-equivalent. Thus,
   $G_{i,\D} = G_{\D}$.

   So, the induction assumption and \eqref{eqn:pionne} imply \eqref{eqn:stab}.

 Let us prove   \eqref{eqn:pione} for the above map $\Phi$. For
 $\g  \in \pi_1(S_t'',o_t)$ and $t \le s$, \eqref{eqn:pione} follows from relation \eqref{eqn:stab} and the induction assumption. If $S_{s+1}$ is marked (i.e., $G_2$ is quasi-Fuchsian), then for  $\g \in \pi_1(S_{s+1}'',o_{s+1})$ \eqref{eqn:pione} follows from the property \eqref{eqn:pione} for \qf groups.

This proves statement A in the
case when every loop $\a_j$ divides $S_0$.

\subsection{Proof of statement A in the
case when some loop $\a_k$ is non-dividing}

The set of all liftings of $\a_k $ divides $\D_0$ into regions. Let $\t \D_0$ be one of these regions, and let $G^*$ be the subgroup of $G$ keeping $\t \D_0$ invariant. Let $ \D_0^*$ be that component of $G^*$ containing  $\t \D_0$.

The following statements are proved in \cite{maskit}.

i. The factor $\t \D_0/ G^* =  S_0'$  is the surface $S_0$ "cut'' along $\a_k$.

ii. The set $\D_0^*$ is simply connected and invariant under $G^*$, and $\D_0^*/G^* = S_0^*$ is the surface $S_0'$ with two disjoined punctured discs $C_1$ and $C_2$ attached along the boundary curves formed by the cut.

iii. $$
\Om(G^*)/G^* = S_0^* + \sum_1^s S_j^+,
$$
$$
\Om(G)/G = S_0 + \sum_1^s S_j^+.
$$
  It is important that all the components of the factors of these two groups, except for the first ones, coincide.

  iv. There exist fundamental sets $D, D^*$ for $G, G^*$ respectively in their discontinuity sets so that

  \begin{equation}\label{eqn:dstar}
  D - \D_0 = D^*- \D_0^*.
  \end{equation}

   Statements iii and iv imply that every component $\Delta$ of $G^*$ different from $\D_0^*$ is at the same time a component of $G$,  and the projections $\pi$ and $\pi^*$ on it coincide. Indeed, each point $x\in\Delta$ is $G^*$-equivalent to a unique point
   $y\in D^*-\D_0^*=D - \D_0\subset\Omega(G)$ by \eqref{eqn:dstar}:  $y=g(x)$, $g\in G^*$. Hence, $\Delta\subset\Omega(G)$. The component in $\Omega(G)$ containing $\Delta$
   coincides with $\Delta$, since $\Omega(G)\subset\Omega(G^*)$, as at the end of the previous subsection.

  Conversely, each non-invariant component $\Delta$ of $G$ is
   $G$-equivalent to a non-invariant component of $G^*$. Indeed, each point $x\in\Delta$ is $G$-equivalent to a unique point $y=g(x)$, $g\in G$, and $y\in D - \Delta_0 $.
   By  \eqref{eqn:dstar}, $y \in D^*\setminus\Delta_0^*$. Therefore, the component
   $g(\Delta)$ of the group $G$ is contained in a component of the group $G^*$, by the above inclusion. The latter ambient component is also a component of the group $G$,
   by the above argument. Finally, $g(\Delta)$ is a component of $G^*$.
 This argument also implies that if the above $\Delta$  is itself a non-invariant component of the group $G^*$, then $g\in G^*$.

 Let us prove that
   $G_\Delta= G^*_\Delta$ for every component $\Delta$ of the group $G^*$ different from
   $\Delta_0^*$. Indeed, two points  $x_1,x_2\in \Delta$ are $G$-equivalent, if and only if the corresponding points $y_j=g_j(x_j)\in D - \D_0=D^*-D_0^*$ coincide.
   But in this case they are $G^*$-equivalent, since  $g_j\in G^*$ by the above discussion, as at the end of the previous subsection.

   The associated complex $K$ of the group $G$ differs from that (denoted $K^*$) of the group
   $G^*$ by the connector $c_k$ corresponding to $\alpha_k$.
   Let $S_\tau$ be the component of the  complex $K$
   adjacent to $c_k$. Let $\Phi^*:S_0^*\to K^*$ be  a map
   satisfying \eqref{eqn:stab}, \eqref{eqn:pione} for the group $G^*$: it exists by the induction hypothesis. Then it automatically satisfies \eqref{eqn:stab}, \eqref{eqn:pione} for the group $G$ for every $\gamma\subset S_t''$,   $S_t$ being a marked component with $t\neq\tau$, by the above statements on coincidence of components of groups $G$, $G^*$ and the corresponding stabilizers.

   Let us pass to construction of a homeomorphism $\Phi$
   satisfying \eqref{eqn:stab}, \eqref{eqn:pione} for the group $G$ by a modification of the homeomorphism
   $\Phi^*$ corresponding to $G^*$. Recall that $S_0'$ is the surface $S_0$ cut along the loop $\alpha_k$, and its boundary contains two closed components, denoted $a_1$,
   $a_2$, each identified with $\alpha_k$.
   Let $\wt S_\tau'$ denote the connected component
   adjacent to $\alpha_k$ (from the both sides) of the complement $S_0'\setminus\sqcup_{\ell\neq k}\alpha_\ell$. The surface $S_\tau'\subset S_0^*$ corresponding to the group $G^*$ is obtained from the surface $\wt S_\tau'$ by pasting
   two disjoint punctured disks $C_i$ to $a_i$, $i=1,2$.
  Let $u_k\Subset U_k \subset S_0$ be  small nested tubular neighborhoods of $\a_k$.
  For every $i=1,2$ let $u_{ik}$ denote the connected component adjacent to $C_i$
  of the complement $u_k\setminus\a_k$. Let
 $U_{ik}$ denote the connected component adjacent to $u_i$
 of the complement $(U_k - u_k)$.  Set
   $$\Phi|_{S_0\setminus U_k}=\Phi^*.$$
  Let us extend thus defined $\Phi$ continuously to $U_k$
  as in the previous subsection
  so that $\Phi|_{U_{ik}}$ be a squeezing homeomorphism of an open annulus
  $U_{ik}$ onto the punctured disk $\wt C_i:=\Phi(U_{ik}\cup a_i\cup C_i)$ agreeing with
   $\Phi$ on $\partial U_k$ and $\Phi(u_k)=c_k$.
Thus defined $\Phi$ satisfies
statement \eqref{eqn:pione} for the group $G$. Indeed, for marked components
$S_t^+$ different from $S_\tau$  \eqref{eqn:pione} holds  by induction hypothesis and construction. For the component $S_\tau$ (if marked) \eqref{eqn:pione} holds again by
the above statements on coincidence of components and stabilizers of the groups
$G$ and $G^*$ and since
 $K$ and $K^*$ differ by a connector only.
  This proves statement A in the
case when some loop $\a_k$ is non-dividing.

Thus the Main Lemma, and together with it, existence of factor and the uniqueness theorem for completely non-degenerate B-groups are proved.

\section{Acknowledgments}
We are grateful to Mikhail Kapovich, Cyril Lecuire and Andrei Vesnin for helpful discussions and remarks.

\end{document}